\documentclass[12pt]{amsart}

\usepackage{amsmath,xspace,amssymb,mathrsfs}
\usepackage{color}

\input xy
\xyoption{all}
\xyoption{2cell}
\UseAllTwocells
\CompileMatrices

\newcommand{\Spec}{\operatorname{Spec}}
\renewcommand{\phi}{\varphi}

\newcommand{\Ker}{\operatorname{Ker}}

\newcommand{\Ima}{\operatorname{Im}}

\newcommand{\Max}{\operatorname{Max}}

\newcommand{\Min}{\operatorname{Min}}

\newtheorem{proposition}{Proposition}[section]
\newtheorem{lemma}[proposition]{Lemma}

\newtheorem{corollary}[proposition]{Corollary}
\newtheorem{theorem}[proposition]{Theorem}

\theoremstyle{definition}

\newtheorem{remark}[proposition]{Remark}

\begin{document}

\title[Flat topology and applications]{New methods toward the patch and flat topologies with applications}

\author[A. Tarizadeh]{Abolfazl Tarizadeh}
\address{Department of Mathematics, Faculty of Basic Sciences, University of Maragheh \\
P. O. Box 55136-553, Maragheh, Iran.
 }
\email{ebulfez1978@gmail.com}

\date{}

\begin{abstract} In this paper, we use elementary and simple ideas which are based on the significant applications of the power set ring to rebuild and study the patch topology on the prime spectrum from a completely different and new point of view. Specially, the proof of a major result in the literature on the comparison of topologies greatly simplified and shortened. Then we develop more natural and simple methods to obtain the flat topology and its properties. In particular, it is shown that the minimal spectrum of a commutative ring is quasi-compact with respect to the flat topology. Then as an application of this result, all of the related results of Kaplansky, Anderson, Gilmer-Heinzer, Bahmanpour-Khojali-Naghipour and Naghipour on the finiteness of the minimal primes are deduced as special cases of this result. Dually, a similar result is also obtained for maximal ideals. \\
\end{abstract}

\maketitle

\section{Introduction}

Historically, the patch and flat topologies were studied by \cite[Chap 3, Ex 27-30]{Atiyah}, \cite{Fontana et al.}, \cite{Fontana 2}, \cite[p. 333-339]{Grothendieck}, \cite{Hochster}, \cite[\S V]{Picavet}, \cite{Ruza-Vielma} and \cite{Tarizadeh}. In the literature, the patch topology is also called the constructible topology, the flat topology is also called the inverse topology. In \cite[Theorem 3.17]{Tarizadeh}, we fixed a mistake of \cite[Remark 2.20 (f)]{Fontana et al.}. Our method in \cite{Tarizadeh} to construct the flat topology is completely different from the Hochster's approach, but it is more similar to the approach of \cite{Fontana et al.}. In \cite{Tarizadeh}, we only used the standard methods of commutative algebra to build the patch and flat topologies. \\

In the present paper in \S3, we introduce a topology on the prime spectrum by significant using of the power set ring and prove that this topology and the patch topology are the same. The ideas which are used here are elementary, more natural, so interesting, and have the potential that can be developed further to apply for other purposes. In particular, we greatly simplified and shortened the approaches of \cite{Fontana 2} to demonstrate that the patch and ultrafilter topologies are the same things, see Theorem \ref{Theorem V}. \\

In \S 4, we adopt a new way to obtain the flat topology and its properties. This method is more natural and simple than the approaches of \cite[\S2]{Fontana et al.}, \cite{Hochster} and \cite[\S3]{Tarizadeh}. In particular, Theorem \ref{Theorem NMTPF I} is obtained. Using this result, then Theorem \ref{prop 677} is deduced which states that ``the minimal spectrum of a commutative ring is quasi-compact with respect to the flat topology''. Then as an application of this result, it is shown that a ring $R$ has finitely many minimal primes if and only if for each $\mathfrak{p}\in\Min(R)$ there exists a finitely generated ideal $I\subseteq\mathfrak{p}$ of $R$ such that $\Min(R)\cap V(I)$ is a finite set, see Theorem \ref{coro 3}. It should be noted that the minimal spectrum of a ring is not necessarily quasi-compact with respect to the Zariski topology. Hence, in the literature, there are various results on the finiteness of the minimal spectrum which most of them have been proved by fairly hard and long methods. All of these results specially \cite{Anderson}, \cite[Theorem 2.1]{Kamal}, \cite[Theorem 1.6]{Gilmer-Heinzer}, \cite[Theorem 88]{Kaplansky} and \cite{Naghipour} are easily deduced as special cases of Theorem \ref{prop 677} or Theorem \ref{coro 3}. Dually, a similar result is also obtained for maximal ideals, see Theorem \ref{theorem 2}.
\\

\section{Preliminaries}

We need the following material in the sequel. \\

In this paper, all rings are commutative. If $\phi:A\rightarrow B$ is a morphism of rings then the induced map $\Spec(B)\rightarrow\Spec(A)$ is denoted by $\phi^{\ast}$ which maps each prime ideal $\mathfrak{p}$ of $B$ into $\phi^{-1}(\mathfrak{p})$. \\

Let $R$ be a ring. A subset $E$ of $\Spec(R)$ is said to be a \emph{patch closed} subset of $\Spec(R)$ if there exists a ring map $\phi:R\rightarrow A$ such that $E=\Ima\phi^{\ast}$. It is well known that the collection of patch closed subsets of $\Spec(R)$, as the closed subsets, forms a topology over $\Spec(R)$, see e.g. \cite[Proposition 2.3]{Tarizadeh}. It is called the \emph{patch topology}. \\

If $f$ is a member of a ring $R$ then $D(f)=\{\mathfrak{p}\in\Spec(R): f\notin\mathfrak{p}\}$ and $V(f)=\Spec(R)\setminus D(f)$. \\

\begin{proposition}\label{Remark I} The collection of $D(f)\cap V(I)$ where $f\in R$ and $I$ runs through the set of finitely generated ideals of $R$ forms a basis for the patch opens of $\Spec(R)$. \\
\end{proposition}

{\bf Proof.} To see its proof consider \cite[Remark 2.5]{Tarizadeh}. $\Box$ \\

If $X$ is a set then its power set $\mathcal{P}(X)$ together with the symmetric difference $A+B=(A\cup B)\setminus (A\cap B)$ as the addition and the intersection $A.B=A\cap B$ as the multiplication form a commutative ring whose zero and unit are respectively the empty set and the whole set $X$. The ring $\mathcal{P}(X)$ is called the \emph{power set ring} of $X$. If $f:X\rightarrow Y$ is a function then the map $\mathcal{P}(f):\mathcal{P}(Y)\rightarrow\mathcal{P}(X)$ defined by $A\rightsquigarrow f^{-1}(A)$ is a morphism of rings. In fact, the assignments $X\rightsquigarrow\mathcal{P}(X)$ and $f\rightsquigarrow\mathcal{P}(f)$ form a faithful contravariant functor from the category of sets to the category of Boolean rings. We call it the \emph{power set functor}. \\

If $J$ is an ideal of $\mathcal{P}(X)$ and $A_{1},...,A_{n}\in J$ then it is easy to see that $\bigcup\limits_{i=1}^{n}A_{i}\in J$. Also, $\mathcal{P}(X)$ is a Boolean ring and so every prime ideal of $\mathcal{P}(X)$ is a maximal ideal. \\

Let $Y\subseteq X$ be two sets, $M$ a maximal ideal of $\mathcal{P}(X)$ and let $N$ be the extension of $M$ under the canonical ring map $\mathcal{P}(X)\rightarrow\mathcal{P}(Y)$ which is induced by the above inclusion. Then $N$ is a maximal ideal of $\mathcal{P}(Y)$ if and only if $Y\notin N$. \\

Let $R$ be a ring. The set of minimal prime ideals of $R$ is denoted by $\Min(R)$. Also, the set of maximal ideals of $R$ is denoted by $\Max(R)$. \\

A ring map $\phi:R\rightarrow A$ is said to be of finite type if $A$ is a finitely generated $R-$algebra via $\phi$. \\

\section{The patch topology from a new perspective}

In this section we study the patch topology from a new perspective.
The ideas of this section were inspired by the works \cite{Fontana 2} and \cite{Ruza-Vielma}. \\

Let $R$ be a ring and let $X$ be a subset of $\Spec(R)$. Then the map $\phi:R\rightarrow\mathcal{P}(X)$ given by $f\rightsquigarrow D(f)\cap X$ is not necessarily a ring map. This map is clearly multiplicative and preserves units, but it is not additive.
Indeed, $D(f)+D(g)\subseteq D(f+g)$. The inclusion may be the strict. For example, in the ring of integers $\mathbb{Z}$ if we take $f=g=1$, then $D(f+g)=D(2)=\Spec(\mathbb{Z})
\setminus\{2\mathbb{Z}\}\neq\emptyset$ but $D(f)+D(g)=\emptyset$. \\

In spite of this, it is interesting to note that if $M$ is a maximal ideal of $\mathcal{P}(X)$ then $M^{\ast}:=\phi^{-1}(M)=\{f\in R: D(f)\cap X\in M\}$ is a prime ideal of $R$ and the map $\eta_{X}:\Spec\mathcal{P}(X)\rightarrow\Spec(R)$ given by $M\rightsquigarrow M^{\ast}$ is continuous since for each $f\in R$,  $\eta_{X}^{-1}\big(D(f)\big)= D(A)$ where $A=D(f)\cap X$. Also notice that if $x:=\mathfrak{p}\in X$ then $\mathfrak{m}^{\ast}_{x}=\mathfrak{p}$ where $\mathfrak{m}_{x}=\mathcal{P}(X\setminus\{x\})$ is a maximal ideal of $\mathcal{P}(X)$. Therefore $X\subseteq\Ima\eta_{X}$. This observation leads us to the following definition. We call the subset $X$ an \emph{ind-power} subset of $\Spec(R)$ if $X=\Ima\eta_{X}$. Thus, $X$ is an ind-power subset of $\Spec(R)$ if and only if $M^{\ast}\in X$ for all $M\in\Spec\mathcal{P}(X)$. \\

\begin{theorem}\label{Theorem IV} The collection of ind-power subsets of $\Spec(R)$, as the closed subsets, forms a topology over $\Spec(R)$. \\
\end{theorem}

{\bf Proof.} Clearly the empty set and the whole set $\Spec(R)$ are the ind-power subsets of $\Spec(R)$. It remains to show that the ind-power subsets of $\Spec(R)$ are stable under finite unions and arbitrary intersections. Let $(X_{i})$ with $i=1,...,n$ be a finite family of ind-power subsets of $\Spec(R)$. Let $M$ be a maximal ideal of $\mathcal{P}(X)$ where $X=\bigcup\limits_{i=1}^{n}X_{i}$. We have to show that $M^{\ast}\in X$. There exists some $k$ such that $X_{k}\notin M$. \\
Then consider the ring map $\mathcal{P}(X)\rightarrow\mathcal{P}(X_{k})$ induced by the inclusion $X_{k}\subseteq X$. Let $M_{k}$ be the extension of $M$ under this ring map which is a maximal ideal of $\mathcal{P}(X_{k})$. To conclude the assertion it suffices to show that $M^{\ast}=M^{\ast}_{k}$. If $f\in M^{\ast}$ then $D(f)\cap X_{k}\in M_{k}$ and so $f\in M^{\ast}_{k}$. Conversely, if $f\in M^{\ast}_{k}$ then there exists some $A\in M$ such that $D(f)\cap X_{k}=A\cap X_{k}$. This yields that $D(f)\cap X_{k}\in M$. If $D(f)\cap X\notin M$ then $X_{k}.(D(f)\cap X)=D(f)\cap X_{k}\notin M$. But this is a contradiction. Thus $f\in M^{\ast}$. Now let $(X_{i})$ be an arbitrary family of ind-power subsets of $\Spec(R)$. Let $M$ be a maximal ideal of $\mathcal{P}(X)$ where $X=\bigcap\limits_{i}X_{i}$. Again as above, we have to show that $M^{\ast}\in X$. Let $M_{i}$ be the contraction of $M$ under the canonical ring map $\mathcal{P}(X_{i})\rightarrow\mathcal{P}(X)$ which is induced by the inclusion $X\subseteq X_{i}$. It is easy to see that $M^{\ast}=M^{\ast}_{i}$ for all $i$. Hence, $M^{\ast}\in X$. $\Box$ \\

We call the resulted topology in Theorem \ref{Theorem IV}, the \emph{ind-power topology}. \\

The proof of the following result is the culmination of this section. \\

\begin{theorem}\label{Theorem V} The patch and ind-power topologies over $\Spec(R)$ are the same. \\
\end{theorem}

{\bf Proof.} First note that the ind-power topology is finer than the patch topology. To see this it suffices to show that for each $f\in R$ then both $D(f)$ and $V(f)$ are the ind-power subsets of $\Spec(R)$, see Proposition \ref{Remark I}. But these are clearly the ind-power subsets. The patch topology is Hausdorff. Hence, to conclude the assertion it suffices to show that the ind-power topology is quasi-compact. To prove this, let $X:=\Spec(R)=\bigcup\limits_{i\in I}U_{i}$ where each $U_{i}$ is an open subset of $\Spec(R)$ with respect to the ind-power topology. We claim that $\Spec\mathcal{P}(X)=\bigcup\limits_{i\in I}D(U_{i})$. Because, if $M$ is a maximal ideal of $\mathcal{P}(X)$ then there exists some $k$ such that $M^{\ast}\in U_{k}$. We show that $M\in D(U_{k})$. If $U_{k}\in M$ then $Y:=X\setminus U_{k}\notin M$ since every element of $\mathcal{P}(X)$ is an idempotent. Let $N$ be the extension of $M$ under the canonical ring map $\mathcal{P}(X)\rightarrow\mathcal{P}(Y)$ which is induced by the inclusion $Y\subseteq X$. It is easy to see that $N$ is a maximal ideal of $\mathcal{P}(Y)$ and $M^{\ast}=N^{\ast}$. It follows that $M^{\ast}\in Y$ since $Y$ is an ind-power subset of $\Spec(R)$. But this is a contradiction. Hence, $M\in D(U_{k})$. This establishes the claim. But we may write  $\Spec\mathcal{P}(X)=\bigcup\limits_{i=1}^{n}D(U_{i})$ since the Zariski topology is quasi-compact. This yields that $\Spec(R)=\bigcup\limits_{i=1}^{n}U_{i}$. $\Box$ \\

In particular, we obtain the following non-trivial result. \\

\begin{corollary}\label{Corollary I} The ind-power subsets of $\Spec(R)$ are precisely the patch closed subsets of $\Spec(R)$. \\
\end{corollary}

{\bf Proof.} It is an immediate consequence of Theorem \ref{Theorem V}. $\Box$ \\

\begin{remark} Let $X$ be a set. Then it is easy to see that the map $M\rightsquigarrow \mathcal{P}(X)\setminus M=\{A\in\mathcal{P}(X): A^{c}\in M\}$ is a homeomorphism from $\Spec\mathcal{P}(X)$ onto $\mathscr{F}(X)$, the space of ultrafilters on $X$ equipped with the Stone topology. Recall that the collection of $d(A)=\{F\in\mathscr{F}(X):A\in F\}$ with $A\in\mathcal{P}(X)$ forms a basis for the opens of the Stone topology over $\mathscr{F}(X)$. As an application of this identification, it can be shown that for a ring $R$ then the ind-power subsets of $\Spec(R)$ are precisely the ultrafilter closed subsets of $\Spec(R)$ in the sense of \cite[Definition 1]{Fontana 2}.  \\
\end{remark}

\section{More natural way toward the flat topology}

In this section, we introduce a new way to build the flat topology. This method is more natural and simple than the approaches of \cite[\S2]{Fontana et al.}, \cite{Hochster} and \cite[\S3]{Tarizadeh}. Then we extract all of the properties of the flat topology using only this method. \\

Let $X$ be a set and let $\mathcal{B}$ be a collection of subsets of
$X$ such that the following two conditions hold. \\
$\mathbf{(i)}$ $X=\bigcup\limits_{B\in\mathcal{B}}B$. \\
$\mathbf{(ii)}$ If $x\in B\cap B'$ for some $B,B'\in\mathcal{B}$, then there is some $B''\in\mathcal{B}$ such that $x\in B''\subseteq B\cap B'$. \\
Then there exists a unique topology over $X$ such that the collection $\mathcal{B}$ forms a basis for its opens. Using this basic fact, then we get the following result. \\

\begin{theorem}\label{Theorem NewI} Let $R$ be a ring. Then there exists a unique topology over $\Spec(R)$
such that the collection of $V(I)$ where $I$ runs through the set of finitely generated ideals of $R$ forms a basis for its opens. \\
\end{theorem}

{\bf Proof.} Clearly $V(0)=\Spec(R)$. If $I$ and $J$ are ideals of $R$ then $V(I)\cap V(J)=V(IJ)$ and $IJ=(x_{i}y_{j}:i=1,...,n; j=1,...,m)$ is a finitely generated ideal of $R$ whenever $I=(x_{1},...,x_{n})$ and $J=(y_{1},...,y_{m})$ are finitely generated ideals. $\Box$ \\

The resulted topology in Theorem \ref{Theorem NewI}, is called the \emph{flat topology}. \\

\begin{theorem}\label{Theorem NMTPF I} The flat topology over $\Spec(R)$ is quasi-compact. \\
\end{theorem}

{\bf Proof.} By Theorem \ref{Theorem NewI}, the flat topology is coarser than the ind-power topology. In the proof of Theorem \ref{Theorem V}, we showed that the ind-power topology is quasi-compact. $\Box$ \\

Recall that a subset $E$ of $\Spec(R)$ is called stable under the generalization if $\mathfrak{p}\in E$ and $\mathfrak{q}$ is a prime ideal of $R$ such that $\mathfrak{q}\subseteq\mathfrak{p}$, then $\mathfrak{q}\in E$. \\

\begin{theorem}\label{Theorem II} Let $E\subseteq\Spec(R)$. Then $E$ is a flat closed subset of $\Spec(R)$ if and only if it is a patch closed subset of $\Spec(R)$ and is stable under the generalization. \\
\end{theorem}

{\bf Proof.} Let $E$ be a flat closed subset of $\Spec(R)$. Clearly it is patch closed. Take $\mathfrak{p}\in E$ and let $\mathfrak{q}$ be a prime ideal of $R$ such that $\mathfrak{q}\subseteq\mathfrak{p}$. If $\mathfrak{q}\notin E$ then by Theorem \ref{Theorem NewI}, there exists a (finitely generated) ideal $I$ of $R$ such that $\mathfrak{q}\in V(I)\subseteq\Spec(R)\setminus E$. But this is a contradiction. Hence, $E$ is stable under the generalization. Conversely, let $E$ be a patch closed subset of $\Spec(R)$ which is stable under the generalization.
There exists a morphism of rings $\phi:R\rightarrow A$ such that $E=\Ima\phi^{\ast}$. We show that $U:=\Spec(R)\setminus E$ is a flat open of $\Spec(R)$. First note that if $\mathfrak{p}$ is a prime ideal of $R$ then it is easy to see that $\mathfrak{p}\in U$ if and only if $\mathfrak{p}A=A$. Thus if $\mathfrak{p}\in U$ then we may write
$1_{A}=\sum\limits_{i=1}^{n}a_{i}\phi(f_{i})$ where $a_{i}\in A$ and $f_{i}\in\mathfrak{p}$ for all $i$. This yields that $\mathfrak{p}\in V(I)\subseteq U$ where $I=(f_{1},...,f_{n})$ is finitely generated ideal of $R$. Thus by Theorem \ref{Theorem NewI}, $U$ is a flat open of $\Spec(R)$. $\Box$ \\

Consider the canonical ring map $\pi:R\rightarrow S^{-1}R$ where $S$ is a multiplicative subset of a ring $R$. Then clearly $\Ima\pi^{\ast}=\{\mathfrak{p}\in\Spec(R): \mathfrak{p}\cap S=\emptyset\}$ is stable under the generalization. Thus by Theorem \ref{Theorem II}, $\Ima\pi^{\ast}$ is a flat closed subset of $\Spec(R)$. Off course this fact holds more generally, see Proposition \ref{Theorem III}. \\

If $\mathfrak{p}$ is a prime ideal of $R$ then the flat closure of this point in $\Spec(R)$, denoted by $\Lambda(\mathfrak{p})$, is induced by the canonical ring map $\pi:R\rightarrow R_{\mathfrak{p}}$. That is, $\Lambda(\mathfrak{p})=\Ima\pi^{\ast}=
\{\mathfrak{q}\in\Spec(R):\mathfrak{q}\subseteq\mathfrak{p}\}$. By contrast, we know that the Zariski closure of this point in $\Spec(R)$ is induced by the canonical ring map $R\rightarrow R/\mathfrak{p}$. \\

\section{Finiteness of minimal and maximal spectra}

The maximal spectrum of a ring is quasi-compact with respect to the Zariski topology. Dually, we have the following result.\\

\begin{theorem}\label{prop 677} The minimal spectrum of a ring $R$ is quasi-compact with respect to the flat topology. \\
\end{theorem}

{\bf Proof.} It suffices to show that every open covering of the minimal spectrum has a finite refinement. Suppose $\Min(R)\subseteq\bigcup\limits_{i}U_{i}$ where $U_{i}$ is a flat open subset of $\Spec(R)$ for all $i$. Let $\mathfrak{p}$ be a prime ideal of $R$. Then there exists a minimal prime ideal $\mathfrak{q}$ of $R$ such that $\mathfrak{q}\subseteq\mathfrak{p}$. We know that the flat closure of $\{\mathfrak{p}\}$ in $\Spec(R)$ is equal to $\{\mathfrak{p}'\in\Spec(R): \mathfrak{p}'\subseteq\mathfrak{p}\}$. Therefore $\mathfrak{q}\in\overline{\{\mathfrak{p}\}}$.
But there exists some $k$ such that $\mathfrak{q}\in U_{k}$.  It follows that $\mathfrak{p}\in U_{k}$. Therefore $\Spec(R)=\bigcup\limits_{i}U_{i}$. But by Theorem \ref{Theorem NMTPF I}, $\Spec(R)$ is quasi-compact with respect to the flat topology. Hence, the assertion is deduced. $\Box$ \\

The following result improves the main results of \cite[Theorem 2.1]{Kamal} and \cite[Theorem on page 232]{Naghipour}. Also, our proof is short and simple. \\

\begin{theorem}\label{coro 3} For a ring $R$ the following statements are equivalent. \\
$\mathbf{(i)}$ $\Min(R)$ is a finite set. \\
$\mathbf{(ii)}$ If $\mathfrak{p}\in\Min(R)$ then there exists a finitely generated ideal $I$ of $R$ such
that $\Min(R)\cap V(I)=\{\mathfrak{p}\}$. \\
$\mathbf{(iii)}$ If $\mathfrak{p}\in\Min(R)$ then there exists a finitely generated ideal $I\subseteq\mathfrak{p}$ of $R$ such that $\Min(R)\cap V(I)$ is a finite set. \\
\end{theorem}

{\bf Proof.} $\mathbf{(i)}\Rightarrow\mathbf{(ii)}:$ If $\Min(R)$ is a finite set then by the prime avoidance theorem, there exists an element $f\in\mathfrak{p}$ such that $f\notin\bigcup
\limits_{\substack{\mathfrak{q}\in\Min(R),\\
\mathfrak{q}\neq\mathfrak{p}}}\mathfrak{q}$.
It follows that $\Min(R)\cap V(f)=\{\mathfrak{p}\}$.
$\mathbf{(ii)}\Rightarrow\mathbf{(iii)}:$ Obvious. \\
$\mathbf{(iii)}\Rightarrow\mathbf{(i)}:$ Using Theorem \ref{prop 677} and hypothesis, then there exist a finite number $I_{1},...,I_{n}$ of (finitely generated) ideals of $R$ such that $\Min(R)\subseteq\bigcup\limits_{k=1}^{n}V(I_{k})$ and each $U_{k}:=\Min(R)\cap V(I_{k})$ is a finite set. Hence, $\Min(R)=\bigcup\limits_{k=1}^{n}U_{k}$ is a finite set. $\Box$ \\

\begin{corollary}\cite[Theorem on page 232]{Naghipour} Let $R$ be a ring with the property that for each $\mathfrak{p}\in\Min(R)$ then there exists a finitely generated ideal $I\subseteq\mathfrak{p}$ of $R$ such that $\Min(R/I)$ is a finite set. Then $\Min(R)$ is a finite set. \\
\end{corollary}

{\bf Proof.} If $\Min(R/I)=\{\mathfrak{p}_{1}/I,...,\mathfrak{p}_{n}/I\}$ then $\Min(R)\cap V(I)\subseteq\{\mathfrak{p}_{1},...,\mathfrak{p}_{n}\}$. Thus by Theorem \ref{coro 3} (iii), $\Min(R)$ is a finite set. $\Box$ \\

\begin{remark} In regarding with the above result, note that if $I$ is an ideal of a ring $R$ then the map $\Min(R)\cap V(I)\rightarrow\Min(R/I)$ given by $\mathfrak{p}\rightsquigarrow\mathfrak{p}/I$ is an injective map. If moreover, $V(I)$ is stable under the generalization (e.g.
$I$ is a pure ideal i.e., the canonical ring map $R\rightarrow R/I$ is a flat ring map) then the above map is also surjective. \\
\end{remark}

\begin{corollary} Let $R$ be a ring. Then $\Min(R)$ is a finite set if and only if $\Min(R)$ with the flat topology is a discrete space. \\
\end{corollary}

{\bf Proof.} It is an immediate consequence of Theorem \ref{coro 3}. $\Box$ \\

The following result easily yields all of the results \cite{Anderson},  \cite[Theorem 1.6]{Gilmer-Heinzer} and \cite[Theorem 88]{Kaplansky}. \\

\begin{corollary}\label{coro 1} Let $R$ be a ring such that every minimal prime is the radical of a finitely generated ideal of $R$. Then $R$ has finitely many minimal primes. \\
\end{corollary}

{\bf Proof.} If $\mathfrak{p}$ is a minimal prime of $R$ then by the hypothesis there exists a finitely generated ideal $I$ of $R$ such that $\mathfrak{p}=\sqrt{I}$. It follows that $\Min(R)\cap V(I)=\Min(R)\cap V(\sqrt{I})=\{\mathfrak{p}\}$. Then apply Theorem \ref{coro 3}. $\Box$ \\

\begin{remark} Let $\{R_{i} : i\in I\}$ be a family of domains. For each $i$, $\pi_{i}^{-1}(0)$ is a minimal prime of $R=\prod\limits_{i\in I}R_{i}$ where $\pi_{i}:R\rightarrow R_{i}$ is the canonical projection. Because, suppose there exists a prime ideal $\mathfrak{p}$ of $R$ such that $\mathfrak{p}\subset\pi^{-1}_{i}(0)$. Pick $a\in\pi^{-1}_{i}(0)\setminus\mathfrak{p}$. Then $ab=0$ where $b=(\delta_{i,j})_{j\in I}$. It follows that $b\in\mathfrak{p}$, a contradiction.
The minimal prime $\pi_{i}^{-1}(0)$ is principal since it is generated by the sequence $(1-\delta_{i,j})_{j\in I}$. Using this and Corollary \ref{coro 1} then we find a minimal prime $\mathfrak{p}$ in every infinite direct product ring $\prod\mathbb{Z}$ such that $\mathfrak{p}\neq\pi_{i}^{-1}(0)$ for all $i$. Moreover, $\mathfrak{p}$ is not the radical of a finitely generated ideal. Note that the set of minimal primes of $\prod\mathbb{Z}$ is in bijective correspondence to the set of all ultrafilters on the index set, see \cite[Proposition 1]{Ronnie et al}. Hence, the number of the minimal primes of an infinite direct product $\prod\mathbb{Z}$ is equal to $2^{2^{\kappa}}$ where $\kappa$ is the cardinality of the index set. \\
\end{remark}

\begin{proposition}\label{Prop 565} Let $R$ be a ring. Then $\Max(R)$ is a finite set if and only if  for every maximal ideal $\mathfrak{m}$ of $R$ there exists some $f\in R\setminus\mathfrak{m}$ such that $\Max(R)\cap D(f)=\{\mathfrak{m}\}$. \\
\end{proposition}

{\bf Proof.} If $\Max(R)$ is a finite set then for each maximal ideal $\mathfrak{m}$ of $R$ there exists some $f\in\bigcap
\limits_{\substack{\mathfrak{m}'\in\Max(R),\\
\mathfrak{m}'\neq\mathfrak{m}}}\mathfrak{m}'$ such that $f\notin\mathfrak{m}$. The converse is deduced from the fact that the maximal spectrum of a ring is quasi-compact with respect to the Zariski topology. $\Box$ \\

\begin{corollary} Let $R$ be a ring. Then $\Max(R)$ is a finite set if and only if $\Max(R)$ with the Zariski topology is a discrete space. \\
\end{corollary}

{\bf Proof.} It is an immediate consequence of Proposition \ref{Prop 565}. $\Box$ \\

\begin{lemma}\label{Lemma 656} Let $S$ be a multiplicative subset of a ring $R$ such that the canonical ring map $R\rightarrow S^{-1}R$ is of finite type. Then there exists some $f\in S$ such that $S^{-1}R=R[1/f]$. \\
\end{lemma}

{\bf Proof.} There exists a finite subset $\{r_{1}/f_{1},...,r_{n}/f_{n}\}$ of $S^{-1}R$ such that each element $r/s\in S^{-1}R$ can be written as $r/s=g(r_{1}/f_{1},...,r_{n}/f_{n})$ for some $g(x_{1},...,x_{n})\in R[x_{1},...,x_{n}]$. It follows that there is a natural number $N$ such that $r/s=r'/f^{N}$ where $f=f_{1}f_{2}...f_{n}$ and $r'\in R$. $\Box$ \\

As a dual of Corollary \ref{coro 1}, we have the following result.\\

\begin{theorem}\label{theorem 2} Let $R$ be a ring with the property that for each maximal ideal $\mathfrak{m}$ the canonical map $\pi:R\rightarrow R_{\mathfrak{m}}$ is injective and of finite type. Then $R$ has finitely many maximal ideals.\\
\end{theorem}

{\bf Proof.} By Lemma \ref{Lemma 656}, there is some $f\in R\setminus\mathfrak{m}$ such that $R_{\mathfrak{m}}=R[1/f]$.
It follows that $\Max(R)\cap D(f)=\{\mathfrak{m}\}$. Because let $\mathfrak{m}'$ be a maximal ideal of $R$ such that $\mathfrak{m}'\neq\mathfrak{m}$. Then we may choose some $g\in\mathfrak{m}'\setminus\mathfrak{m}$. In the ring $R_{\mathfrak{m}}$ we may write $1/g=r/f^{n}$ for some $n\geq1$. Thus there exists some $h\in R\setminus\mathfrak{m}$ such that $h(f^{n}-rg)=0$. It follows that $f^{n}=rg$ because $\pi$ is injective.
Thus $f\in\mathfrak{m}'$. Therefore by Proposition \ref{Prop 565}, $R$ has finitely many maximal ideals. $\Box$ \\

We conclude this section by proving the following result which shows that the injectivity hypothesis of Theorem \ref{theorem 2} is more natural to happen. \\

\begin{proposition} Let $R$ be a ring. Then for each maximal ideal $\mathfrak{m}$ of $R$ the canonical map $\pi:R\rightarrow R_{\mathfrak{m}}$ is injective if and only if the set of zero-divisors of $R$ is contained in the Jacobson radical of $R$. \\
\end{proposition}

{\bf Proof.} ``$\Rightarrow$'': Suppose there exists a zero-divisor $f\in R$ such that $f\notin\mathfrak{J}(R)$. Thus there exists a maximal ideal $\mathfrak{m}$ of $R$ such that $f\notin\mathfrak{m}$. But there exists some non-zero $g\in R$ such that $fg=0$. So $g\in\Ker\pi$. But this is a contradiction since $\Ker\pi=\{0\}$. Conversely, let $\mathfrak{m}$ be a maximal ideal of $R$ and $h\in\Ker\pi$. Thus there exists some $h'\in R\setminus\mathfrak{m}$ such that $hh'=0$. If $h\neq0$ then $h'$ is a zero-divisor of $R$ and so $h'\in\mathfrak{J}(R)\subseteq\mathfrak{m}$. This is a contradiction and we win. $\Box$ \\

\section{More on the flat topology}

The first proof of the following technical lemma, uses the Gabriel localization theory. Gabriel's approach in Bourbaki \cite[Chap 2. p. 158, Les exercices 16 \`{a} 25]{Bourbaki} is a very nice introduction to the Gabriel localization theory. Also, we refer the interested reader to \cite{Akiba}, \cite{Nastasescu-Popescu}, \cite{Popescu-Spircu} and \cite{Stenstrom}. \\

\begin{lemma}\label{Lemma I} Let $\phi:R\rightarrow A$ be a morphism of rings such that $E:=\Ima\phi^{\ast}$ is stable under the generalization. Then there exists a flat ring map $\psi:R\rightarrow B$ such that $E=\Ima\psi^{\ast}$. \\
\end{lemma}

{\bf First proof.} Let $\mathscr{F}$ be the set of ideals $I$ of $R$ such that $IA=A$. Then $\mathscr{F}$ is an idempotent topologizing system (Gabriel topology) over $R$. Let $R_{\mathscr{F}}$ be the Gabriel localization of $R$ with respect to the system $\mathscr{F}$ which is a commutative $R-$algebra whose the structure morphism is denoted by $\pi:R\rightarrow R_{\mathscr{F}}$. It can be shown that $\pi$ is a flat ring map, and $\phi$ factors uniquely through it i.e., there exists a unique morphism of rings $\psi:R_{\mathscr{F}}\rightarrow A$ such that $\phi=\psi\circ\pi$. This in particular yields that $E\subseteq\Ima\pi^{\ast}$. Conversely, if $\mathfrak{p}\in\Ima \pi^{\ast}$ then $\mathfrak{p}\notin\mathscr{F}$ and so $\mathfrak{p}A\neq A$. Thus there exists a prime ideal $\mathfrak{q}$ of $A$ such that $\mathfrak{p}\subseteq\phi^{-1}(\mathfrak{q})\in E$. Hence, $\mathfrak{p}\in E$ since $E$ is stable under the generalization. Therefore $E=\Ima\pi^{\ast}$. \\
{\bf Second proof.} By Theorem \ref{Theorem II}, $E$ is a flat closed subset of $\Spec(R)$. By a similar argument as applied in the construction of the patch topology, see \cite[Proposition 2.3]{Tarizadeh}, we can show that the collection of $\Ima\theta^{\ast}$ where $\theta:R\rightarrow S$ is a ``flat'' ring map, as the closed subsets, forms a topology over $\Spec(R)$. Then it can be shown that this topology is precisely the flat topology, see \cite[Theorem 3.2]{Tarizadeh}. Hence, there exists a flat ring map $\psi:R\rightarrow B$ such that $E=\Ima\psi^{\ast}$. $\Box$ \\

\begin{proposition}\label{Theorem III} Let $R$ be a ring. Then the flat closed subsets of $\Spec(R)$ are precisely of the form $\Ima\phi^{\ast}$ where $\phi:R\rightarrow A$ is a flat ring map. \\
\end{proposition}

{\bf Proof.} If $\phi:R\rightarrow A$ is a flat ring map then $\Ima\phi^{\ast}$ is stable under the generalization because it is well known that every flat ring map has the going-down property, see e.g. \cite[Theorem 9.5]{Matsumura}. Thus by Theorem \ref{Theorem II}, $\Ima\phi^{\ast}$ is flat closed of $\Spec(R)$. The second proof of Lemma \ref{Lemma I}, also shows that $\Ima\phi^{\ast}$ is flat closed of $\Spec(R)$. For the reverse implication see Lemma \ref{Lemma I}.  $\Box$ \\

In practice, Theorem \ref{Theorem II} is much more efficient than Proposition \ref{Theorem III}. But the cost that we pay for Proposition \ref{Theorem III}, has the worth that it reveals a new aspect of the flat topology, it also justifies the word ``flat'' in terminology of the flat topology. \\

\end{document}